\documentclass[11pt]{elsart}
\reversemarginpar

\usepackage{amsmath} 
\usepackage{amsfonts}
\usepackage{amssymb}
\usepackage{graphicx}
\def\bbbr{{\rm I\!R}}
\newtheorem{theorem}{Theorem}
\newtheorem{problem}{Problem}
\newtheorem{remark}{Remark}
\newtheorem{lemma}{Lemma}
\newcommand{\proof}{{\em Proof. }}
\def\qed{\rightline{$\square$}}
\newtheorem{example}{Example}
\def\BB{{\mathcal B}}
\def\IBB{{\mathcal B}^{-1}}
\def\ITBB{{\mathcal B}^{-T}}

{\bf}{}

\begin{document}
\begin{frontmatter}

\title{Robust control of uncertain multi-inventory systems via Linear Matrix Inequality}

\author[bauso]{D. Bauso,}
\author[Giarre]{L. Giarr\'e}
\author[pesenti]{and R. Pesenti\thanksref{Some}}
\address[bauso]{ Dipartimento di Ingegneria Informatica\\
       Universit\`a di Palermo, Viale Delle Scienze, 90128 Palermo,   Italy.   \\
    {\em dario.bauso@unipa.it}}
\address[Giarre]{   Dipartimento di Ingegneria dell'Automazione e dei Sistemi  \\
        Universit\`a di Palermo, Viale Delle Scienze, 90128 Palermo,   Italy.   \\
    {\em giarre@unipa.it}}
\address[pesenti]{Dipartimento di Matematica Applicata -
Universit\`a ``Ca' Foscari'' di Venezia, Dorsoduro 3825/e, 30123
Venezia, Italy.\\
{\em pesenti@unive.it}}

\thanks[Some]{Corresponding author R. Pesenti, Research supported by PRIN ``Advanced control and
identification techniques for innovative applications'', and PRIN
``Analysis, optimization, and coordination of logistic and
production systems''.}

\small{
\begin{abstract}
We consider a continuous time linear multi--inventory system with
unknown demands bounded within ellipsoids and controls bounded
within ellipsoids or polytopes. We address the problem of
$\varepsilon$-stabilizing the inventory since this implies some
reduction of the inventory costs. The main results are certain
conditions under which $\varepsilon$-stabilizability is possible
through a saturated linear state feedback control. All the results
are based on a Linear Matrix Inequalities (LMIs) approach and on
some recent techniques for the modeling and analysis of polytopic
systems with saturations. 
\end{abstract}}
\end{frontmatter}


\section{Introduction}We
consider a continuous time linear multi--inventory system with
unknown demands bounded within ellipsoids and controls bounded
within ellipsoids or polytopes. The system is  modelled as a first
order one integrating the discrepancy between controls and demands
at different sites (buffers). Thus, the state represents the
buffer levels. We wish to study conditions under which the state
can be driven within an a-priori chosen target set through a
saturated linear state feedback control. Let $\varepsilon$ be a
maximal dimension of the target set, the above problem corresponds
to $\varepsilon$-stabilizing the state.

Motivations for $\varepsilon$-stabilizing the state derive from the
benefits associated to keeping the state and consequently also the
inventory costs bounded. This work is in line with some recent
literature on robust optimization~\cite{AP06,BT06} and
control~\cite{BBP06} of inventory systems. Here as well as
in~\cite{BBP06} we focus on saturated linear state feedback
controls since such controls arise naturally in any system with
bounded controls.

The main results of this work can be summarized as follows.
Initially we introduce the necessary and sufficient conditions for
the $\varepsilon$-stabilizability in the form of an inclusion between
convex sets. In the case where both demands and controls are
bounded within polytopes, it is well known that verifying such
conditions is NP-hard~\cite{MC96}. Here, we prove that
verification becomes easy when both demands and controls are
bounded within ellipsoids (we will refer to it as the
\emph{ellipsoidal case}). This is possible by rewriting the
inclusion between ellipsoids in terms of
unconstrained quadratic maximization. 

For the ellipsoidal case, we first characterize invariant sets
through a fourth degree condition. As verifying such a condition
is difficult, we then propose the best quadratic approximation of
the same condition. We proceed by describing the region of
linearity of the control and conclude by providing LMI conditions
on the target set under which the saturated control
$\varepsilon$-stabilizes the system. The case where demands are
bounded within ellipsoids and controls are bounded within
polytopes (we will refer to it as the \emph{polytopic case}) is an
open problem and we propose certain sufficient LMI conditions to
solve it.


All the results are based on a Linear Matrix Inequalities (LMIs)
approach in line with the recent work~\cite{B06} on
inventory/manufacturing systems. In particular, when addressing
the politopic case, we use the same technique provided
in~\cite{GT01} to rewrite the model with saturations in polytopic
form. Once we do this, we can apply the LMI analysis covered in
the book~\cite{BEFB} for polytopic systems.

This paper is arranged as follows. In Section~\ref{sec:pf}, we
formulate the problem. In Section~\ref{sec:nsc}, we introduce
necessary and sufficient conditions for the admissibility of the
problem. In Sections~\ref{sec:ec} and~\ref{sec:pc} we study the
problem with ellipsoidal and polytopic constraints respectively.
Finally, in Section~\ref{sec:c}, we draw some
conclusions.


\section{Problem Formulation}\label{sec:pf}
Consider the continuous time linear multi--inventory system
\begin{equation}
\dot x(t) = B u(t) - w(t), \label{system}
\end{equation}
where $x(t) \in \bbbr^n$ is a vector whose components are the
buffer levels, $u(t)  \in \bbbr^m$ is the controlled flow vector,
$B \in \mathbb{Q}^{n\times m}$, with $m \geq n$ and $rank(B) = n$
is the controlled process matrix and $w(t) \in \bbbr^n$ is the
unknown demand. To model backlog $x(t)$ may be less than zero.
Demands are bounded within ellipsoids, i.e., 
\begin{equation}
\label{eq:d}w(t)  \in  {\mathcal W}=\{w \in \mathbb R^n: w^T R_w
w\leq 1\}.\end{equation}
In a first case, in the following
referred as \emph{ellipsoidal case}, controls are bounded within
ellipsoids,
\begin{equation}\label{eq:u}u(t)\in
{\mathcal U}=\{u \in \mathbb R^m: u^T R_u u \leq 1\}.
\end{equation}
In a second case, in the following referred as \emph{polytopic
case}, controls are bounded within polytopes
\begin{equation}\label{eq:u1}u(t) \in 
{\mathcal U}=\{u \in \mathbb R^m:\, u^- \leq u \leq u^+\}
\end{equation}
with assigned $u^+$, $u^-$.

For any positive definite matrix $P\in \mathbb R^{n \times n}$,
define the function $V(x)=x^TPx$ and the ellipsoidal target set
$\Pi=\{x\in \bbbr^n: V(x)\leq 1\}$. In addition, for any matrix
$K\in \mathbb R^{n \times n}$, define as saturated linear state
feedback control any policy
\begin{equation}\label{eq:SatPolicy}
u=-sat\{Kx\}=\left\{
\begin{array}{ll}
    -Kx & \mbox{if } Kx \in \mathcal{U}\\
u(x) \in \partial \mathcal{U} & \mbox{otherwise}\end{array}
\right.\end{equation} where hereafter $\partial F$ indicates the frontier of a
given set $F$.

\begin{problem}\label{prob:1}($\varepsilon$-stabilizing)
Consider a system~(\ref{system}) in the ellipsoidal or polytopic
case. Find conditions on the positive definite matrix $P\in
\mathbb R^{n \times n}$, under which there exists a saturated linear state feedback control $u=-sat\{Kx\}$
such that
it is possible to drive the
state $x(t)$ within the target set $\Pi$.
\end{problem} Solving the above problem corresponds to
$\varepsilon$-stabilizing the state $x$ within $\Pi$.

\begin{example}\label{sec:IllEx} Throughout this paper
we consider, as illustrative example, the graph with one node and
two arcs depicted in Fig.~\ref{Graph1}. The incidence matrix is
$B=[1 \quad 1]$.
\begin{figure}\label{Graph1}
\begin{center}
\includegraphics[width=5cm]{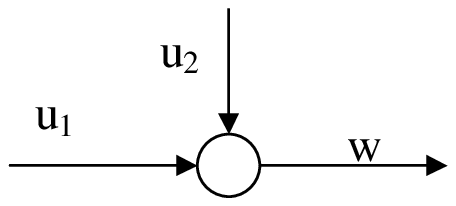}
\caption{Graph with one node and two arcs.}
\end{center}
\end{figure}
The continuous time dynamics is
$$
\dot x (t) =\underbrace{[1\quad 1]}_{B}
\underbrace{\left[\begin{array}{cc}u_1(t)\\u_2(t)\end{array}\right]}_u-
w= u_1(t) + u_2(t) - w(t),
$$
with demand bounded in the ellipsoid $$w^2\leq 1$$ and with the
following either ellipsoidal or polytopic constraints on the
control~$u$
\begin{equation}\label{eq:exelly}
(u_1+u_2)^2\leq 1,  
\end{equation}
\begin{equation}
\label{eq:expoly}-2\leq u_1\leq 3,\,  -2\leq u_2\leq 1. 
\end{equation}
Finally, the target set is the sphere of
unitary radius $\Pi = \{x \in \mathbb R: x^2 \leq 1\}$.
\end{example}


\section{Stability necessary and sufficient conditions}\label{sec:nsc}

System~(\ref{system}) is $\varepsilon$-stabilizable if and only if
for all $w \in {\mathcal W}$, there exists $u \in int\{{\mathcal
U}\}$ such that $Bu = w$ (see, e.g., \cite{BRU97}). For the short
of notation, the previous condition is usually expressed as
\begin{equation}B{\mathcal U} \supset {\mathcal W}.\label{eq:BUW}\end{equation}
Deciding whether~(\ref{eq:BUW}) holds is NP-hard, when ${\mathcal
U}$ and ${\mathcal W}$ are polytopes. Here, we prove that
verifying~(\ref{eq:BUW}) becomes easy when both ${\mathcal U}$ and
${\mathcal W}$ are ellipsoids. Observe that we can rewrite $Bu =
w$ as $u_\BB = \IBB w-\IBB N u_N$, where $B = [\BB | N]$ being
$\BB$ a basis of $B$ and $N$ the remaining columns of $B$,
correspondingly $u_\BB$ are the $n$ components of $u$ associated
to the basis $\BB$ and $u_N$ are  the $m-n$ components of $u$
associated to the columns in $N$.

As we observe that~(\ref{eq:BUW}) is
equivalent to
$$\max_{w \in {\mathcal W}} \min_{u \in \mathbb{R}^{m}: Bu = w} uR_u u< 1,$$
Condition~(\ref{eq:BUW}) holds if and only if
\begin{eqnarray}
\max_{w \in {\mathcal W}} \min_{u_N \in \mathbb{R}^{m-n}}&&
f(u_\BB(w,u_N),u_N) = \nonumber \\&&= \left[w^T \ITBB - u_N^T
N^T\ITBB | u_N^T\right]R_u \left[\begin{array}{c} \IBB w-  \IBB N
u_N\\ u_N\end{array}\right]< 1\label{eq:minmax1}\end{eqnarray}

 When we consider the illustrative example in Section
\ref{sec:IllEx}, we have $\BB = [1]$, $N = [1]$ then problem
(\ref{eq:minmax1}) becomes
\begin{eqnarray}\max_{-1 \leq w\leq 1} \min_{u_2 \in \mathbb{R}}&& f(u_\BB(w,u_2),u_2) = \nonumber \\&&= \left[w -u_2 | u_2\right]\left[\begin{array}{cc} 1 & 0 \\ 0 & 1\end{array}\right] \left[\begin{array}{c} w-u_2\\ u_2\end{array}\right]= (w-u_2)^2+u_2^2< 1\label{eq:minmax1es}\end{eqnarray}

Now consider, function $ f(u_\BB(w,u_N),u_N)$. It is a differentiable convex
function in $u_N$. Then, for any $w \in {\mathcal W}$ we can
analytically determine the best response $u_N^*(d) = arg\min_{u_N
\in \mathbb{R}^{m-n}} f(u_\BB(w,u_N),u_N)$, by imposing
$$
\nabla_{u_N}f(u_\BB(w,u_N),u_N)=2\left[-N^T\ITBB
| I\right]R_u \left[\begin{array}{c} \IBB w - \IBB N u_N\\
u_N\end{array}\right] = 0,$$  where $I$ is the $(m-n)\times (m-n)$
identical matrix. We obtain$$u_N^*(w) =
-\underbrace{\left(\left[-N^T\ITBB |
I\right]R_u \left[\begin{array}{c} - \IBB N\\
I\end{array}\right]\right)^{-1} \left[\begin{array}{c} -
N^T\ITBB | I\end{array}\right]R_u\left[\begin{array}{c} \IBB\\
0\end{array}\right]}_{M}w=-Mw,$$ where $0$ is the $(m-n)\times n$
null matrix. In the example under consideration, we have $$u_2^*(w)
= -\left(\left[-1 |1\right]\left[\begin{array}{cc} 1 & 0 \\ 0 &
1\end{array}\right] \left[\begin{array}{c} -1\\
1\end{array}\right]\right)^{-1}\left[-1
|1\right]\left[\begin{array}{cc} 1 & 0 \\ 0 & 1\end{array}\right]
\left[\begin{array}{c} 1\\ 0\end{array}\right]w=\frac{w}{2}.$$

For any $w \in {\mathcal W}$ the minimal value of $
f(u_\BB(w,u_N),u_N)$ is $$f(u_\BB(w,u_N^*(w)),u_N^*(w))=
w^{*T}\Phi w^*,$$ where  \begin{equation}
\Phi=\underbrace{[\ITBB +
M^TN^T\ITBB | -M^T]}_{H^T}R_u \underbrace{\left[\begin{array}{c} \IBB + \IBB NM\\
-M\end{array}\right]}_{H}=H^T R_u H
\label{eq:DefH}
\end{equation}  is a positive definite $n
\times n$ matrix, as $M$ is full rank. So far, we have shown that
we can find the optimal value of problem~(\ref{eq:minmax1}) by
solving problem
\begin{equation}\max_{w \in \mathcal W}= w^{T}\Phi w,\label{eq:minmax2}\end{equation}
and checking that the optimal value is less than one.

We are ready to
observe that problem~(\ref{eq:minmax2}) is easy as it reduces to
determining the eigenvectors of an $n \times n$ matrix.

\begin{theorem}\label{th:InclusionTheorem}
System~(\ref{system}) is $\varepsilon$-stabilizable if and only if
$w^{*T}\Phi w^* < 1$, for all $w^*$ eigenvectors associated to the
maximum eigenvalue of matrix $R_w^{-1}\Phi $ whose weighted
quadratic norm $w^{*T} R_w w^*$ is equal to 1.
\end{theorem}
\proof As  $w^{T}\Phi w$ is convex, its optimal value $w^*$ lays on the
frontier $\partial {\mathcal W}$ of the set ${\mathcal W}$, i.e.,
for $w^{*T} R_w w ^* = 1$.  Imposing the Karush Kuhn Tucker first
order optimality condition, we obtain $2(\Phi-\lambda R_w)w^* =0$.
Then the optimal values of $w^*$ are some of the matrix
$R_w^{-1}\Phi $ eigenvectors whose weighted quadratic norm $w^{*T}
R_w w^*$ is equal to 1. In particular, $w^*$ are the eigenvectors
associated to the maximal eigenvalues of $R_w^{-1}\Phi $.

\qed

In the example under consideration $\Phi =
\left[\frac{1}{2}\right]$ and $w^*= \pm 1$ then $w^{*T}\Phi w^* =
\frac{1}{2}< 1$, hence the associated system is
$\varepsilon$-stabilizable.

In the following we discuss for which initial state  the system is
certainly $\varepsilon$-stabilizable through a (pure) linear state
feedback control; hence we show that if we saturated the previous
linear policy the system is
$\varepsilon$-stabilizable for any  initial state. 

\section{Ellipsoidal constraints}\label{sec:ec}

Let us start by considering only the constraints~(\ref{eq:d}) on
$w$ and neglect the ellipsoidal constraints~(\ref{eq:u}) on $u$.
Among the saturated linear state feedback
control~(\ref{eq:SatPolicy}) we prove that we can solve
Problem~\ref{prob:1} using controls of type $u=sat\{-kHx\}$, with
$k \in \mathbb{R}$ and $H \in \mathbb{R}^n$ as defined
in~(\ref{eq:DefH}). Note that matrix $H$ is a right inverse of
$B$, that is $BH = I$. We motivate the choice of $u=-sat\{kHx\}$
with $H$ as defined in~(\ref{eq:DefH}) as such a control describes
the best response of $u$ under the worst $w$ as proved in the
previous section. Also, note that the scalar $k \in \mathbb{R}$
must be lower than a certain value, which means that we cannot use
a bang-bang control. This is motivated by the following reason. If
we use a control $u=sat\{-kHx\}$, then the necessary and
sufficient condition (\ref{eq:BUW}) becomes
\begin{equation}\label{eq:blin}B{\mathcal U}_{lin} \supset {\mathcal W}\end{equation}
where
$${\mathcal U}_{lin}=\{u \in \mathbb R^m: \,u=-kHx,\,k^2x^TH^TR_uHx\leq 1\}.$$
Following the derivation of (\ref{eq:minmax2}) in the previous
Section, we have that (\ref{eq:blin}) holds if and only if $$k^2
w^{*T}\Phi w^*<1.$$ For $k=1$ the above condition holds true as it
reduces to (\ref{eq:minmax2}). Obviously, the value $\hat k
=\sqrt{\frac{1}{w^{*T}\Phi w^*}}$ is an upper bound for $k$,
namely, we must choose $k$ such that $k< \hat k$ if we wish the
necessary and sufficient condition~(\ref{eq:blin}) be satisfied.

With the above considerations in mind, we can conclude that the
dimensions of the target $\Pi$ where it is possible to drive the
state are lower bounded.

Denote by $\lambda_{max}(Z)$ the maximum eigenvalue of a given
matrix $Z$. In the following theorem we prove that $\dot V(x) < 0$
within a given set (\emph{invariant} set). This result will allow
exploiting  $V(x)$ as a Lyapunov function to prove the convergence
to the target set $\Pi$.

\begin{theorem}\label{eq:nouconst0}
Consider system~(\ref{system}) subject to the only ellipsoidal
constraints~(\ref{eq:d}) on $w$, and controlled via linear state
feedback $u=-kHx$, with $H$ such that $BH = I$. Then condition
$\dot V<0$ holds if and only if
\begin{equation}\label{eq:quad}k^2(x^TPx)^2 -x^TPR_w^{-1}Px > 0.\end{equation}
\end{theorem}

\proof For $H$ such that $BH = I$, condition $\dot V<0$ is
equivalent to
\begin{equation}\label{eq:vdotnewrp} 2kx^TPx+2w^TPx > 0.\end{equation}
We aim at proving that $\dot V<0$ holds for any $x$ external to an
appropriate smooth closed surface. To do this, we look for an $x
\in \mathbb{R}^{n}$ inducing a solution strictly greater than zero
for the following problem
\begin{equation}\label{eq:vdotext} \min_{w \in \mathcal{W}} \zeta(x,w) = 2kx^TPx+2w^TPx.\end{equation}
As $\zeta(x,w) $ is linear in $w$, the optimal $w^*$ must lay on
the boundary of set $\mathcal{W}$. The Karush Kuhn Tucker
conditions impose that $Px = -\lambda R_w w^*$ for some $\lambda
\geq 0$, that is $w^* =
-\frac{1}{\lambda}R_w^{-1}Px$. Note that being $P$ full rank, it necessarily holts that$\lambda \not = 0$ for all $x  \not = 0$.  Then, $
\zeta(x,w^*) = 2kx^TPx -\frac{2}{\lambda}x^TPR_w^{-1}Px > 0$. As
$w^*$ lays on the boundary of $\mathcal{W}$, we have $w^{*T}R_w
w^*=\frac{x^TPR_w^{-1}Px}{\lambda^2}=1$ from which
$\lambda=\sqrt{x^TPR_w^{-1}Px}$. Hence, $\zeta(x,w^*) > 0$, and
therefore also~(\ref{eq:vdotnewrp}) holds, if and only if
(\ref{eq:quad}) holds.

\qed

We now exploit  $V(x)=x^TPx$ as a Lyapunov function to prove the convergence
to the target set $\Pi$. We determine under which conditions on $P$ and $k$ we have that $\dot V < 0$ or, equivalently,
 inequality (\ref{eq:quad}) hold for any $x \not \in \Pi$.

When $P=\nu R_w$, (\ref{eq:quad}) becomes
$k^2x^TPx > \nu$. Then, in this case,  we can use $V(x)$  to prove the convergence of the system to $\Pi$ for $k^2 \geq \nu$.

In the following, we consider the general case when  $P\not =\nu R_w$.

\begin{lemma}\label{eq:unit}
Consider system~(\ref{system}) subject to the only ellipsoidal
constraints~(\ref{eq:d}) on~$w$, and controlled via linear state
feedback $u=-kHx$, with $H$ such that $BH = I$. Then, $k^2(x^TPx)^2 -x^TPR_w^{-1}Px > 0$ holds for  any $x \not \in \Pi$
if and only if  $k^2 -x^TPR_w^{-1}Px \geq 0$ holds for any $x \in \partial \Pi$.
\end{lemma}

\proof
(\emph{Necessity}). Assume that there exists $\hat x \in \partial \Pi$ such that $k^2 -x^TPR_w^{-1}Px < 0$. Then, there also exists a ball $Ball(\hat x, r)$ centered in $\hat x$ with a sufficiently small radius $r >0$ such that for all $x \in Ball(\hat x, r)$ we have $k^2 -x^TPR_w^{-1}Px < 0$. This implies that there exist
$x \not \in \Pi$ for which condition (\ref{eq:quad}) does not hold.

(\emph{Sufficiency}). Assume that $k^2 -x^TPR_w^{-1}Px \geq 0$ holds for any $x \in \partial \Pi$. By contradiction, consider $\hat x \not \in \Pi$, i.e., $\hat x^TP\hat x = \rho > 1$, such that $k^2(\hat x^TP\hat x)^2 -\hat x^TPR_w^{-1}P\hat x < 0$, that is $k^2 \rho^2 -\hat x^TPR_w^{-1}P\hat x < 0$. Then, there exists $\tilde x = \frac{\hat x}{\sqrt{\rho}} \in \partial \Pi$ such that $k^2 \rho^2 - \rho \tilde x^TPR_w^{-1}P\tilde x < 0$, that is $k^2 \rho -  \tilde x^TPR_w^{-1}P\tilde x < 0$. This latter result is contradictory as we cannot have $k^2 \rho < \tilde x^TPR_w^{-1}P\tilde x \leq  k^2$, for $\rho > 1$.

\qed


\begin{lemma}\label{eq:nouconst}
Consider system~(\ref{system}) subject to the only ellipsoidal
constraints~(\ref{eq:d}) on~$w$, and controlled via linear state
feedback $u=-kHx$, with $H$ such that $BH = I$. We can use $V(x)$  to prove the convergence of the system to $\Pi$ for $k^2 \geq \lambda_{max}(R_w^{-1}P)$.
\end{lemma}

\proof Condition $k^2 -x^TPR_w^{-1}Px \geq 0$ holds for any $x \in \partial \Pi$ if and only if $\min_{x \in \partial \Pi}\{k^2 -x^TPR_w^{-1}Px \} \geq 0$.
Imposing the Karush Kuhn Tucker first
order optimality condition, we obtain $2(PR_w^{-1}P-\lambda P)x^* =0$.
Then the optimal values of $x^*$ are some of the matrix
$R_w^{-1}P $ eigenvectors whose weighted quadratic norm $x^{*T}
P x^*$ is equal to 1. In particular, $x^*$ are the eigenvectors
associated to the maximal eigenvalues of $R_w^{-1}P$. For vectors $x^*$, condition $k^2 -x^{*T}PR_w^{-1}Px^* \geq 0$  becomes $k^2 -\lambda_{max}(R_w^{-1}P)x^{*T}Px^* \geq 0$, that is $k^2 -\lambda_{max}(R_w^{-1}P) \geq 0$.

\qed

%
%


Observe that the system converges to the target set $\Pi_R=\{x:
k^2 x^TR_wx \leq 1\}$ as any feasible target set $\Pi=\{x: x^TPx
\leq 1\}$, with $k^2 \geq \lambda_{max}(R_w^{-1}P)$ includes
$\Pi_R$. Indeed,  $\Pi \supseteq \Pi_R$ if $x^TPx -  k^2 x^TR_wx =
x^T(P-k^2R_w)x \leq 0$ or  equivalently if $P-k^2R_w \preceq 0$.
In turn, the latter condition  is equivalent to  $R^{-1}_wP-k^2I
\preceq 0$ that certainly holds as $k^2\geq
\lambda_{max}(R_w^{-1}P)$

In the next theorem we introduce the constraints on controls
(\ref{eq:u}). To this end, we need to define the family of
ellipsoid $\Sigma_0(\xi) =\{x \in \mathbb R^n: x^T P x\leq x(0)^T
P x(0):=\xi\}$ parametrized in $\xi \geq 1$.

\begin{theorem}\label{th:ell-lin}
Given system~(\ref{system}) in the ellipsoidal case, we can drive
the state $x(t)$ from any initial value $x(0)\in \Sigma_0(\xi)$ to the
target set $\Pi$ via linear state feedback $u=-kHx$ if the
following conditions hold
\begin{equation}\label{eq:KBthm1new}
k^2  \geq  \lambda_{max}(R_w^{-1}P)\\
\end{equation}
\begin{equation}
\label{eq:KBthm1new2} k^2 \xi\lambda_{max}( P^{-1} \Phi)\leq  1. 
\end{equation}
\end{theorem}

\proof By Lemma~\ref{eq:nouconst}, under condition
(\ref{eq:KBthm1new}) it holds $\dot V(t)< 0$ for all $x(t) \not
\in \Pi$ and then $V(x)$ can be considered as a Lyapunov function
for the convergence of the state to the set  $\Pi$ when the linear
control $u=-kHx$ is implemented. Condition $\dot V(t)< 0$  also
implies that $\Sigma_0(\xi)$ is invariant with respect to the same
linear feedback as $\xi \geq 1$ which means $\Sigma_0(\xi) \supseteq
\Pi$. Then
$$\max_{t\geq 0} u^T(t) R_u u(t) \leq \max_{x\in \Sigma_0(\xi)} k^2 x^T
H^TR_u H x=\max_{x\in \Sigma_0(\xi)} k^2 x^T \Phi
x=k^2\xi\lambda_{max}( P^{-1} \Phi).$$ Therefore the constraint
$u=-kHx(t)\in {\mathcal  U}$ for all $t\geq 0$ is enforced
if~(\ref{eq:KBthm1new2}) holds true.

\qed

%

The following theorem provides a solution to Problem~\ref{prob:1}.
Let us denote by $X$ the set of states $x$ where we can define a
linear control $u(x) = - kHx$, i.e.,  $X = \{x: -kHx \in \mathcal
U\}$. Consider the saturated linear state feedback control of type
\begin{equation}\label{eq:slsf}u(x) = \left\{ \begin{array}{ll}
   -k Hx & \mbox{if } x \in X \\
   -\frac{Hx}{\sqrt{x^TH^TR_uHx}} & \mbox{if } x \not \in X
\end{array} \right. .\end{equation}

\begin{theorem}\label{th:sp2} Consider a system~(\ref{system}) in the ellipsoidal case.
For any positive definite matrix $P\in \mathbb R^{n \times n}$
satisfying condition~(\ref{eq:KBthm1new}), the saturated linear
state feedback control~(\ref{eq:slsf}) drives the state $x(t)$
within the target set $\Pi$ for any initial state
$x(0)$.\end{theorem} \proof By construction, $u(x)$ is a
continuous function with $\mathcal  U$ as codmain. When we use
such a control, we know that $\dot V(x) < 0$ also holds for any $x
\not \in \Pi$, if $\Pi \subset X$ and $k^2 \geq
\lambda_{max}(R^{-1}P)$ (see Lemma~\ref{eq:nouconst}).

First observe that, for all $x \in \partial X$, we have
$x^TPx>k^2x^TH^TR_uHx=1$, where the latter inequality holds as
$\Pi \subset X$. Then, for any $x \not \in  X$, that is for
$k^2x^TH^TR_uHx>1$, we have $ \frac{x^TPx}{x^TH^TR_uHx}>k^2\geq
\lambda_{max}(R^{-1}P)$  since both $x^TPx$ and $x^TH^TR_uHx$ are
positive definite quadratic forms.

In Lemma \ref{eq:nouconst}, we have  proved  that $\dot V(x) < 0$
for $x \in X \setminus \Pi$. Now, we consider $x \not \in X$. We
have $\dot V(x) < 0$ if and only if $-x^TPBu(x)+x^TPw >0$, for all
$w \in \mathcal  W$, that is

\begin{equation}
   \min_{w \in \mathcal  W}\left\{\frac{x^TPx}{\sqrt{x^TH^TR_uHx}}
+x^TPw\right\} >0
   \label{eq:dotVgeneral}
\end{equation}
must hold.  Applying the Karush-Kuhn-Tuker conditions, we
transform (\ref{eq:dotVgeneral}) in
$\frac{x^TPx}{\sqrt{x^TH^TR_uHx}} - \sqrt{x^TP^TR^{-1}_wPx} > 0$.
In turn, the latter inequality holds if $\frac{x^TPx}{x^TH^TR_uHx}
- \lambda_{max}(R^{-1}P) > 0$, as $x^TP^TR^{-1}_wPx \leq
\lambda_{max}(R^{-1}P) x^TPx$. We then conclude that $\dot V(x)<
0$ since $\frac{x^TPx}{x^TH^TR_uHx}> k^2 \geq
\lambda_{max}(R^{-1}P)$.

\qed

 Observe that the saturated linear
state feedback control (\ref{eq:slsf}) is not decentralized in the
sense that the generic $i$th control $u_i$ in general depends on
the demand at different nodes and on the other controls $u_j,$
$j\not = i$. This is due to either the structure of matrix $H$ or
the ellipsoidal constraints~(\ref{eq:u}).

\begin{remark} Consider the two equivalent matrix inequalities on
$P$ and $Q=P^{-1}$,
\begin{equation}\label{eq:lmiQ}(2k-1)P-P R_w^{-1}P  \geq
0,  \quad (2k-1)Q- R_w^{-1}  \leq
0.\end{equation} 
 Trivially, any $P$ satisfying condition
(\ref{eq:KBthm1new}) also satisfies the two above matrix
inequalities.\end{remark} Matrix inequalities of the above form
will be used in the following sections.

\begin{example}Consider the graph depicted in
Fig.~\ref{Graph1}, with one node and two arcs and incidence matrix
$B=[1 \quad 1]$. Controls are subject to ellipsoidal constraints
(\ref{eq:exelly}). Then we have, $R_w=1$, $R_u=I$ and
$\Phi=\frac{1}{2}$. We can stabilize the system within $\Pi = \{x \in \mathbb R: x^2 \leq 1\}$ for any initial state $x(0)\leq \sqrt{2}$ via
a pure linear state feedback $u=-kHx$. To see this take $Q=I$, and
observe that the matrix inequality on $Q$~(\ref{eq:lmiQ}) is
satisfied for any $k\geq 1$. Furthermore, if we assume $k=1$, then
from~(\ref{eq:KBthm1new2}) we must have $k^2=1\leq
\frac{2}{\xi^2}=\frac{2}{x(0)^2}$.
\end{example}


\section{Polytopic constraints}\label{sec:pc} Controls $u$ are
subject to the polytopic constraints~(\ref{eq:u1}). Again, we
study under which conditions we can solve  Problem~\ref{prob:1}
using controls of type $u=-sat\{kHx\}$, with $k \in \mathbb{R}$
and $H \in \mathbb{R}^n$ such that $BH=I$. In this case, we
interpret the $sat\{.\}$ operator as componentwise. More
specifically, we choose the control
\begin{equation}\label{eq:newsat}u_i=sat_{[u_i^-,u_i^+]}\{-kH_{i\bullet}x\},\end{equation} with $H$ such that
$BH = I$, $H_{i\bullet}$ denoting the $i$th row of $H$ and where,
for any given scalar $a$ and $b$
$$
sat_{[a,b]}\{\zeta\} = \left \{
\begin{array} {ll}
b,~~\mbox{if}~~\zeta > b,\\
\zeta, ~~\mbox{if}~~a \leq \zeta \leq b,\\
a,~~\mbox{if}~~\zeta < a.
\end{array}
\right .
$$Henceforth we
omit the indices of the $sat$ function.

Under the control $u=sat\{-kHx\}$, the closed loop dynamics becomes
\begin{equation}\label{eq:dynxpoly}\dot x=B sat\{-kHx\}-w.\end{equation}
Our idea is to rewrite the above dynamics in the following
polytopic form
\begin{equation}\label{eq:polysystem} \dot x=A(t)x(t)-w(t),  \quad   w(t)^T
R_w w(t) \leq 1,\end{equation} where  the time varying matrices
$A(t)$ are expressed as convex combinations of $2^m$ matrices
$A_j$, $j=1,\ldots,2^m$. More precisely the expressions for $A(t)$
are
\begin{equation}\label{eq:degree}A(t)=\sum_{j=1}^{2^m} \sigma_j(t) A_j, \quad
\sum_{j=1}^{2^m}
\sigma_j(t)=1.
\end{equation}

The procedure to compute matrices $A_j$'s is borrowed
from~\cite{GT01} and recalled below. Let us rewrite the control
policy as
$$u_i=sat\{-kH_{i\bullet} x \}=\theta_i(x)(-kH_{i\bullet} x),$$ where
$\theta_i(x)$ are the ``degree of saturation'' of the control
components defined as follows
\begin{equation}\label{eq:theta}
\theta_i(x)=\left\{\begin{array}{cl}\frac{u_i^-}{-kH_{i\bullet}x} & \mbox{if $-kH_{i\bullet}x<u_i^-$}\\
1 & \mbox{if $u_i^-\leq -kH_{i\bullet} x\leq u_i^+$}\\
\frac{u_i^+}{-kH_{i\bullet}x} & \mbox{if
$-kH_{i\bullet}x>u^+$}\end{array}\right..\end{equation} Let
$\underline \theta=[\underline \theta_1,\ldots,\underline
\theta_m]$ be a vector whose components $\underline \theta_i$ are
such that $0\leq \underline \theta_i\leq 1$ and represent lower
bounds of $\theta_i(x(t))$, for $t\geq 0$. Lower bounds depend on $x(0)$ and can be
computed as $\underline \theta_i= \min_{x \in \Sigma_0(\xi)} \theta_i(x)$
where we remind the definition of $\Sigma_0(\xi)=\{x\in \mathbb R^n:
x^TPx\leq x(0)^TPx(0):=\xi\}$. Also define the vector
$\psi^\theta=[\psi_1^\theta,\ldots,\psi_m^\theta]$ with
$\psi_i^\theta=\frac{1}{\underline\theta_i}$ and the associated
portion of the state space $$S(\psi^\theta)=\{x \in \mathbb R^n:
\, -\psi^\theta \leq -kH x \leq \psi^\theta\}.$$ According to the
above definition of the $\theta_i$s we derive that $S(\psi^\theta)
\supseteq \Sigma_0(\xi)$. Note that we can affirm that $\underline \theta_i$ are lower bounds because the state trajectory
never exits $S(\psi^\theta)$ as we will show in the proof of Theorem~\ref{th:th-x}.

Consider now the $2^m$ vectors $\gamma_j\in \{1,\underline
\theta_1\}\times \ldots \times \{1,\underline \theta_m\}$, with
$j=1,\ldots,2^m$. In other words, $\gamma_j$ is an $m$ component
vector with $i$th component $\gamma_{ji}$ taking value $1$ or
$\underline \theta_i$. Then, each matrix $A_j$ can be expressed as
$A_j=-Bkdiag(\gamma_j)H$. Roughly speaking each vector $\gamma_j$
stores the minimum and or maximum degree of saturation of all
control components. Also, note that matrices $A_j$s induce a
partition of $S(\psi^\theta)$ into regions $X_j$, with
$j=1\ldots,2^m$. Each region is defined as the set of state values
such that the control components are saturated with degree of
saturation equal to $\gamma_{ji}$, namely
$$X_j=\{x\in \mathbb
R^n:\,\theta_i(x)=\gamma_{ji},\,i=1,\ldots,m\}.$$ We remind here
that $\gamma_{ji}$ is the $i$th component of $\gamma_j$.

To complete the derivation of the polytopic form
(\ref{eq:polysystem}) it is left to be noted that given any $x(t)
\in S(\psi^\theta)$ we can compute the associated degree of
saturation from~(\ref{eq:theta}) and derive the weights
$\sigma_j(t)$ of the convex combination~(\ref{eq:degree}). All the
results in the rest of this section try to give an answer to
Problem~\ref{prob:1} with respect to the polytopic
system~(\ref{eq:polysystem}). For each $A_j$, let us define a
matrix
$$M_{j}=QA_j^T+A_jQ+\alpha Q+\frac{1}{\alpha}R_w^{-1}
$$ for a given positive and arbitrarily chosen scalar $\alpha$ and let $(\lambda_{j}^r,v_{j}^r)$ with $r\in
\{1,\ldots,n\}$ be the negative eigenvalues and corresponding
eigenvectors of $M_{j}$.
\begin{theorem}\label{th:th-x} Consider
system~(\ref{system}) in the polytopic case. The saturated linear
state feedback control~(\ref{eq:newsat}) drives the state $x(t)$
within the target set $\Pi$ if
\begin{equation}\label{eq:th-x}
X_j\subseteq Span\{v_{j}^r\}, \quad \mbox{for all
$j=1,\ldots,2^n$.}\end{equation}
\end{theorem}
\proof First of all, note that if~(\ref{eq:th-x}) holds true then
$\Sigma_0(\xi)$ is invariant. Consequently,  as $\Sigma_0(\xi)
\subseteq S(\psi^\theta)$ and by definition $x(0) \in
\Sigma_0(\xi)$, we also have that the state trajectory $x(t)$ will
never exit $S(\psi^\theta)$. Now, we must show that $\dot V(x)<0$
for all $x$ and $w$ such that $x \not \in \Pi$, $u\in {\mathcal
U}$ and $w\in {\mathcal  W}$. In formulas, we must have
\begin{equation}\label{eq:vdot}\begin{array}{lll}\dot V(x)&=&\dot x^TPx + x^T
P \dot x = [A(t)x-w]^T  P x + x^T P [A(t)x-
w]~=\\&=&x^TA(t)^TPx+x^TPA(t)x-w^T Px-x^TP w<
0\end{array}\end{equation} for all $x$ and $w$ satisfying
\begin{equation}\label{eq:LMI1-1}1-x^TPx\leq 0\end{equation}
\begin{equation}
\label{eq:LMI1-3}w^TR_ww-1\leq 0.\end{equation}

Using the ${\mathcal  S}$-procedure, we can say that
condition~(\ref{eq:vdot}) is implied by conditions
(\ref{eq:LMI1-1})-(\ref{eq:LMI1-3}) if there exist
$\alpha,\beta\geq 0$, such that for all $x$ and $w$
\begin{equation}\label{eq:LMI2}\left[\begin{array}{c}x \\  w\end{array}\right]^T
\left[\begin{array}{cc}A(t)^TP+PA(t)^T + \alpha P & -P
\\-P &
-\beta R_w\end{array} \right] \left[\begin{array}{c}x  \\
w\end{array}\right]- \alpha + \beta \leq 0.\end{equation}
Trivially it must hold $\beta\leq \alpha$. Assume without loss of
generality $\beta= \alpha$. Remind that $\alpha$ and $\beta$ can
be chosen arbitrarily. After pre and post-multiplying by
$Q=P^{-1}$, the above condition becomes
\begin{equation}\label{eq:LMI2-Q}\left[\begin{array}{c}x \\  w\end{array}\right]^T
\left[\begin{array}{cc}QA(t)^T+A(t)^T Q + \alpha Q & -I
\\-I &
-\alpha R_w\end{array} \right] \left[\begin{array}{c}x  \\
w\end{array}\right]\leq 0.\end{equation}

Now, as the state never leaves the region $S(\psi^\theta)$, i.e.,
$x(t)\in S(\psi^\theta)$, we can always express $A(t)$ as convex
combination of the $A_j$s as in~(\ref{eq:degree}).

By convexity, the above condition is true if it holds, for all
$j=1,\ldots,2^n$,
\begin{equation}\label{eq:LMI2-Qv}
\left[\begin{array}{c}x_{(j)} \\
w_{(s)}\end{array}\right]^T\left[\begin{array}{cc}QA_j^T+A_j^T Q +
\alpha Q & -I
\\-I &
-\alpha R_w\end{array} \right] \left[\begin{array}{c}x_{(j)}  \\
w_{(s)}\end{array}\right]\leq 0.\end{equation}

Using the Shur complement the condition~(\ref{eq:LMI2-Qv}) is
implied by~(\ref{eq:th-x}).

\qed

Stronger conditions are established in the following theorem which
also highlights the dependence of $M_{j}$ on the scalar $\alpha$.
\begin{theorem}\label{th:KBthm2}Consider system~(\ref{system}) in the polytopic
case. The saturated linear state feedback control
(\ref{eq:newsat}) drives the state $x(t)$ within the target set
$\Pi$ if there exists a scalar $\alpha\geq 0$ such that
\begin{equation}\label{eq:KBthm2}M_{j}<
0,\quad \mbox{for all $j=1,\ldots,2^n$}.
\end{equation}
\end{theorem}
\proof Trivially, if we observe that~(\ref{eq:KBthm2}) implies
(\ref{eq:th-x}).

\qed

Both (\ref{eq:th-x}) and (\ref{eq:KBthm2}) are  sufficient, but not necessary, conditions. When
they hold, we are sure that the system state converge to a state strictly included in the target
set $\Pi$. We discuss more on
this topic in the next section.


\subsection{Approximation error}
We wish to estimate the difference in terms of volumes between the
target set $\Pi$ and the target set obtained from
conditions~(\ref{eq:KBthm2}) and we will call such a difference as
\emph{approximation error}. On this purpose, denote by $Q_{j}$ the
matrix of the smallest (in volume) ellipsoid satisfying $M_{j}<0$,
which is given by \begin{equation}\label{eq:app}Q_{j}=\arg \inf_Q
\min_\alpha \{det(Q),\,M_{j}=QA_j^T+A_jQ+\alpha
Q+\frac{1}{\alpha}R_w^{-1}<0\}.\end{equation} To do this, let
matrix $\underline A$ be the matrix $A_j$ with $j=1,\ldots,2^m$
obtained when no controls are saturated and note that the dynamics
associated to this single matrix is the same as if we assumed the
controls unbounded. To be more precise, $\underline A=-BkH$ as all
components of $\gamma_j$ are equal to one. Remind that $\gamma_j$
stores the degree of saturation of each control component. Also
let us define $\underline Q$ the solution of~(\ref{eq:app}) for
$A_j=\underline A$. We do this, as the
target set $\Pi$ within which we can stabilize
the state, must inscribe the ellipsoid defined by $\underline Q$,
i.e.,
$$\Pi \supset \{x\in \mathbb R^n:\,x^T \underline Q^{-1} x\leq 1\}.$$
Similarly, let matrix $\overline A$ be the matrix $A_j$ with
$j=1,\ldots,2^m$ obtained when all controls are saturated at their
lowest degree of saturation. To be more precise, $\overline
A=-Bkdiag([\underline\theta_1,\ldots,\underline \theta_m])H$ as
all components of $\gamma_j$ are equal to $\underline \theta_i$
for $i=1,\ldots,m$. If we also define $\overline Q$ the solution
of~(\ref{eq:app}) for $A_j=\overline A$, the target set $\Pi$ must
be inscribed in the ellipsoid defined by $\overline Q$, namely,
$$\Pi \subset \{x\in \mathbb R^n:\,x^T \overline Q^{-1} x\leq 1\}.$$

The approximation error can be measured by the ratio
$$e=\frac{det(\overline
Q^{-1}) - det(\underline Q^{-1})}{det(\underline Q^{-1})}.$$

\begin{example}Consider the graph depicted in
Fig.~\ref{Graph1}, with one node and two arcs, incidence matrix
$B=[1 \quad 1]$, and target set $\Pi = \{x \in \mathbb R: x^2 \leq 1\}$. Controls are subject to polytopic constraints
(\ref{eq:expoly}). Take $H=[\frac{1}{2}\quad \frac{1}{2}]^T$ and
$k=1$. Then according to~(\ref{eq:theta}) we have (here $x$ is a
scalar)
$$\theta_1(x)=\left\{\begin{array}{cl}\frac{2}{x/2} & \mbox{if
$x/2
> 2$}\\ 1 & \mbox{if $-3 \leq x/2 \leq 2$}\\ -\frac{3}{x/2} & \mbox{if
$x/2 < -3$}
\end{array}\right. \quad \theta_2(x)=\left\{\begin{array}{cl}\frac{2}{x/2} & \mbox{if $x/2 >
2$}\\ 1 & \mbox{if $-1 \leq x/2 \leq 2$}\\ -\frac{1}{x/2} &
\mbox{if $x/2 < -1$}\end{array}\right..$$ If we consider initial
states $x(0)$ satisfying $-10\leq x(0) \leq 10$, possible lower
bounds for the $\theta$'s are $\underline \theta_1=\frac{2}{5}$
and $\underline \theta_2=\frac{1}{5}$. Note that
$S(\psi^{\theta})=\{x \in \mathbb R ^n: -10 \leq x \leq 10\}$.
Vectors $\gamma$'s and matrices $A$'s turn
out to be
\begin{equation}\begin{array}{llll}\gamma_1=[1\quad 1]^T & \gamma_2=[0.4\quad 1]^T
& \gamma_3=[1\quad 0.2]^T & \gamma_4=[0.4\quad 0.2]^T\\
A_1=-2 & A_2=-1.4 & A_3=-1.2 & A_4=-0.6\end{array}.\end{equation}
Dynamics~(\ref{eq:polysystem}) is then
\begin{equation}\dot x = [-\sigma_1(t) 2 - \sigma_2(t) 1.4-
\sigma_3(t) 1.2 - \sigma_4(t) 0.6] x +  w,
\end{equation}
with $\sum_{j=1}^4\sigma_j(t)=1$. Furthermore, we have
\begin{equation*}\begin{array}{ll}M_{1}=[-4+\alpha]Q+\frac{1}{\alpha} &
M_{2}=[-2.8+\alpha]Q+\frac{1}{\alpha}\\
M_{3}=[-2.4+\alpha]Q+\frac{1}{\alpha}&
M_{4}=[-1.2+\alpha]Q+\frac{1}{\alpha}
\end{array}.\end{equation*}
To apply Theorem~\ref{th:th-x} and~\ref{th:KBthm2},
note that
$\underline A=A_4$ and that $M_{4}<0$ implies
consequently $M_{j}<0$ for all $j$. The solution of (\ref{eq:app}), for $j=4$ is $Q_4={\underline{Q}}=\frac{1}{0.36}$ and $\alpha=0.6$, then the approximation error is
$e = \frac{1-0.36}{0.36}=1.78$.

\end{example}

\section{Conclusions and future works}\label{sec:c}
We have addressed the problem of $\varepsilon$-stabilizing the
inventory of a continuous time linear multi--inventory system with
unknown demands bounded within ellipsoids and controls bounded
within ellipsoids or polytopes. Motivations are due to the cost
reduction associated with a bounded inventory. As main results we
have provided certain LMIs conditions under which
$\varepsilon$-stabilizability is possible through a saturated linear
state feedback control. We have also exploited some recent
techniques for the modeling and analysis of polytopic systems with
saturations. 

This work is a continuation of~\cite{BBP06} and is in line with
some recent applications of LMI techniques to
inventory/manufacturing systems~\cite{B06}. In a future work, we
will study the validity in probability of the LMI conditions
derived in this paper. This is in accordance with some recent
literature on \emph{chance LMI constraints} developed in the area
of robust optimization~\cite{BN02,CC05}.


\begin{thebibliography}{1}

\bibitem{AP06}
E. Adida, and G. Perakis, ``A Robust Optimization Approach to
Dynamic Pricing and Inventory Control with no Backorders'',
\emph{Mathematical Programming}, Ser. B vol. 107, 2006, pp.
97--129.

\bibitem{BBP06} D. Bauso, F. Blanchini, R. Pesenti, ``Robust control policies
for multi-inventory systems with average flow constraints'',
\emph{Automatica}, Special Issue on Optimal Control Applications
to Management Sciences, vol. 42, no. 8, pp. 1255-1266, Aug. 2006.

\bibitem{BRU97} F.~Blanchini, F.~Rinaldi and W.~Ukovich,
``A network design problem for a distribution system with
uncertain demands'', \emph{SIAM Journal on Optimization}, 7
(1997), pp.~560--578.

\bibitem{BN02}
A. Ben-Tal, A. Nemirovsky, ``On tractable approximations of
uncertain linear matrix inequalities affected by interval
uncertainty'', \emph{SIAM Journal on Optimization}, vol. 12,
pp.~811--833, 2002.

\bibitem{BT06}
D. Bertsimas, A. Thiele, ``A Robust Optimization Approach to
Inventory Theory '', \emph{Operations Research}, vol. 54, no. 1,
Jan-Feb 2006, pp. 150--168.

\bibitem{B06} E. K. Boukas, ``Manufacturing Systems: LMI Approach'',
\emph{IEEE Transactions on Automatic Control}, vol.~51, no.~6,
pp.~1014-1018, June 2006.

\bibitem{BEFB}
S. Boyd, L. El Ghaoui, E. Feron, and V. Balakrishnan, \emph{Linear
Matrix Inequalities in System and Control Theory}, volume 15 of
\emph{Studies in Applied Mathematics}, Society for Industrial and
Applied Mathematics (SIAM), Philadelphia, PA, 1994.

\bibitem{CC05}
G. Calafiore, M. C. Campi, ``Uncertain Convex Programs: Randomized
Solutions and Confidence Levels'', \emph{Mathematical
Programming}, 102, pp.~25--46, 2005.

\bibitem{GT01} J. M. Gomes da Silva, Jr. and S. Tarbouriech
``Local Stabilization of Discrete-Time Linear Systems with
Saturating Controls: An LMI-based Approach'', \emph{IEEE Trans. on
Automatic Control}, vol.~46, no.~1, pp.~119-124, Jan. 2001.

\bibitem{MC96} S.T.~McCormick,  ``Submodular containment
is hard, even for networks'',   \emph{Operations Research
Letters}, 19 (1996), pp.~95--99.

\bibitem{PG01} J. A. Primbs, M. Giannelli, ``Kuhn-Tucker-Based Stability Conditions for Systems With Saturation'',
\emph{IEEE Transactions on Automatic Control}, vol.~46, no.~10,
pp.~1643-1647, Oct. 2001.

\end{thebibliography}
\end{document}